\catcode`\@=11
\pretolerance=500    \tolerance=1000  \brokenpenalty=5000

\catcode`\:=\active
\def:{\relax\ifhmode\ifdim\lastskip>\z@\unskip\fi\penalty\@M\ \fi\string:}

\catcode`\!=\active
\def!{\relax\ifhmode\ifdim\lastskip>\z@\unskip\fi\kern.2em\fi\string!}

\catcode`\?=\active
\def?{\relax\ifhmode\ifdim\lastskip>\z@\unskip\fi\kern.2em\fi\string?}

\def\^#1{\if#1i{\accent"5E\i}\else{\accent"5E #1}\fi}
\def\"#1{\if#1i{\accent"7F\i}\else{\accent"7F #1}\fi}

\def\Medbreak{\vskip-\lastskip\medbreak}
  
\let\bigf@nt=\tenrm     \let\smallf@nt=\sevenrm

\def\pc#1{\bigf@nt#1\smallf@nt}         \def\pd#1 {{\pc#1} }

\long\def\th#1 #2\enonce#3\endth{%
   \Medbreak
   {\pc#1} {#2\unskip} {\sl #3}\medskip}

\long\def\tha#1 #2\enonce#3\endth{%
   \Medbreak
   {\pc#1} {#2\unskip}\par\nobreak{\sl #3}\medskip}

\frenchspacing
\catcode`\@=12

\overfullrule=0pt
\parindent=0cm

\input xy
\xyoption{all}

\def\dm{\mathop{{\bf DM}^{-}}\nolimits}
\def\dme{\mathop{{\bf DM}^{-,\rm eff}}\nolimits}

\def\dg{\mathop{{\bf DM}_{\rm gm}}\nolimits}
\def\sh{\mathop{\rm Shv}\nolimits}
\def\sp{\mathop{\rm Spec}\nolimits}
\def\ho{\mathop{\rm Hom}\nolimits}
\def\alt{\mathop{\rm Alt}\nolimits}

\def\om{\mathop{\bf \Omega}\nolimits}
\def\smc{\mathop{\rm Smcor}\nolimits}
\def\sm{\mathop{\rm Sm}\nolimits}

\def\ke{\mathop{\rm Ker}\nolimits}
\def\Gal{\mathop{\rm Gal}\nolimits}

\def\M{{\bf M}}

\def\limind{\mathop{\oalign{\rm lim\cr
\hidewidth$\longrightarrow$\hidewidth\cr}}}

\def\build#1_#2^#3{\mathrel{\mathop{\kern 0pt#1}\limits_{#2}^{#3}}}

\def\hfl#1#2{\smash{\mathop{\hbox to 12mm {\rightarrowfill}}
\limits^{\scriptstyle#1}_{\scriptstyle#2}}}

\def\dfleche#1#2{\smash{\mathop{\hbox to
9mm{\rightarrowfill}}\limits^{\scriptstyle#1}_{\scriptstyle#2}}}

 \def\gfleche#1#2{\smash{\mathop{\hbox to
9mm{\leftarrowfill}}\limits^{\scriptstyle#1}_{\scriptstyle#2}}}

\def\relrightarrow{\mathrel{\hbox to 9mm{\rightarrowfill}}}
\def\llhooksymb{\lhook\joinrel\relrightarrow}
\def\llhook{\mathrel{\llhooksymb}}

\newtoks\ref
\newtoks\AUTHOR
\newtoks\TITLE
\newtoks\BOOKTITLE
\newtoks\PUBLISHER
\newtoks\YEAR
\newtoks\JOURNAL
\newtoks\VOLUME
\newtoks\PAGES
\newtoks\NUMBER
\newtoks\SERIES
\newtoks\ADDRESS

\def\pointir{\unskip . --- \ignorespaces}

\def\book{\leavevmode
{[\the\ref]\enspace}%
\the\AUTHOR\pointir
\the\TITLE, 
{\sl\the\SERIES},
{\bf \the\VOLUME},
\the\PUBLISHER,
\the\ADDRESS,
({\the\YEAR}).
\smallskip
\filbreak}

\def\article{\leavevmode
{[\the\ref]\enspace}%
\the\AUTHOR\pointir
\the\TITLE,
{\sl\the\JOURNAL},
{\bf\the\VOLUME},
({\the\YEAR}),
{no \the\NUMBER},
\the\PAGES.
\smallskip
\filbreak}

\def\articles{\leavevmode
{[\the\ref]\enspace}%
\the\AUTHOR\pointir
\the\TITLE,
{\sl\the\JOURNAL},
{\bf\the\VOLUME},
({\the\YEAR}),
\the\PAGES.
\smallskip
\filbreak}

\def\prep{\leavevmode
{[\the\ref]\enspace}%
\the\AUTHOR\pointir
\the\TITLE,
{\sl en pr\'eparation}.
\smallskip
\filbreak}

\def\incollection{\leavevmode
{[\the\ref]\enspace}%
\the\AUTHOR \pointir
\the\TITLE,
{\bf in }\the\BOOKTITLE, 
{\sl \the\SERIES },
{\bf\the\VOLUME} ,
\the\PUBLISHER ,
\the\ADDRESS, 
({\the\YEAR}), 
\the\PAGES.
\smallskip
\filbreak}

\def\biblio#1{\vglue 15mm\centerline{\bf#1}\vskip 10mm}
\let\+=\tabalign
\def\signature#1\endsignature{\vskip 15mm minus 5mm\rightline{\vtop{#1}}}

\centerline {\bf LE COMPLEXE MOTIVIQUE DE DE RHAM}
\bigskip\bigskip

\centerline {Florence LECOMTE et Nathalie WACH}
\bigskip\bigskip
\medskip \centerline {Strasbourg}
\vskip0.7cm	
{\leftskip 10pt \rightskip 10 pt {\sl
ABSTRACT. Sur un corps de caract\'eristique $0$, nous construisons un complexe 
motivique de De Rham qui permet de g\'en\'eraliser la cohomologie de De Rham 
d'une vari\'et\'e lisse \`a tout motif au sens de Voevodsky.}\par}

\bigskip
{\bf Introduction}

\smallskip
Sur un corps de caract\'eristique $0$, nous disposons gr\^ace aux travaux de 
Voevodsky de la cat\'egorie triangul\'ee tensorielle $\dm (k)$ des 
complexes motiviques et de sa sous-cat\'egorie pleine $\dg (k)$ des 
motifs g\'eom\'etriques.
La cat\'egorie $\dg (k)$ est engendr\'ee 
comme cat\'egorie triangul\'ee par les motifs des sch\'emas 
projectifs lisses. Utilisant cette description, Annette Huber [H00] a d\'efini un foncteur de $\dg (k)$ vers la cat\'egorie
des r\'ealisations mixtes qu'elle a construite dans [H95].

Dans ce travail, nous proposons une d\'efinition de la r\'ealisation de De Rham d'un complexe motivique  $\M$ de $\dm (k)$ qui  \'etend la composante de De Rham des r\'ealisations de A.Huber. Plus pr\'ecis\'ement, nous construisons un complexe motivique $\om^\bullet$, ind-objet de la cat\'egorie $\dm (k)$,
appel\'e complexe motivique de De Rham et qui, par troncature par
tout entier $n$,  permet d'obtenir une suite d'objets $\tau_{\leq n} \om^\bullet$ de $\dm (k)$. On d\'efinit alors
la r\'ealisation de De Rham d'un objet  $\M$ de $\dm (k)$ par sa $p$-i\`eme composante, en posant,
$$H^p_{DR} (\M ) = \limind_n \ho_{\dm(k)} (\M , \tau_{\leq n} \om^\bullet [p]),$$
pour tout entier $p$. Le complexe $\om^\bullet$ est d\'efini par
$$\matrix{
{\bf \Omega} ^p(X) &= 0\hfill &{\rm si}\  p\leq 0, \cr
&=\Gamma (X,{\cal O}_X) &{\rm si}\  p=0,     \cr
&=\Gamma (X,\Omega^p_X) &{\rm si}\  p\geq 1, \cr
}$$
 pour $X$ sch\'ema lisse sur $k$ et $\Omega^p_X$ le faisceau des 
$p$-formes diff\'erentielles sur $X$.
Il repr\'esente la cohomologie de De Rham dans le sens que pour $X$ sch\'ema lisse sur $k$ de dimension $\leq n$, on a
$${\bf H}_{DR}^p(X)=\ho_{\dme(k)}({\bf M}(X),\tau_{\leq n}{\bf \Omega}^\bullet[p]),$$
o\`u ${\bf M}(X)$ est le motif associ\'e \`a $X$.

\smallskip
Les complexes motiviques 
sont construits \`a  partir des faisceaux de Nisnevich avec transferts;
la difficult\'e essentielle de ce travail consiste \`a munir de transferts les pr\'efaisceaux $X \mapsto \Omega^p_X (X)$ des 
$p$-formes diff\'erentielles des sch\'emas lisses sur $k$. 
\tha THEOREME
\enonce Soient $X$ et $Y$ deux sch\'emas lisses de type fini sur k et soit $W$
une correspondance finie au sens de Voevodsky [V00] de $X$ vers $Y$.
Alors, pour $p\geq1$, il existe un morphisme 
$$\Omega^p_Y (Y) \longrightarrow \Omega^p_X (X)$$
compatible aux compositions des correspondances.
\endth

La d\'emonstration repose principalement sur les travaux de Suslin et Voevodsky 
[SV96]. Mais contrairement \`a ({\it loc. cit.}) les faisceaux $\Omega^p$ 
consid\'er\'es ici ne v\'erifient pas de descente 
galoisienne. C'est pourquoi nous avons recours aux travaux de C. Knighten sur
 les diff\'erentielles de Zariski [K73].

\smallskip
Cet article constitue une premi\`ere \'etape dans notre travail de 
g\'en\'eralisation des r\'ealisations des motifs g\'eom\'etriques aux 
complexes motiviques [LW].

\medskip
{\bf 1   Un complexe de faisceaux de Nisnevich}
\smallskip
Pour chaque entier $p\geq1$, on consid\`ere le pr\'efaisceau sur la cat\'egorie des sch\'emas lisses sur $k$ d\'efini par
$$\matrix{ {\bf \Omega}^p: &\sm(k)&\longrightarrow &{\cal A}b\cr
 &X &\longmapsto &\Omega^p_X(X)\quad.}$$

\th PROPOSITION 1.1 
\enonce  ${\bf \Omega}^p$ est un faisceau pour la topologie de Nisnevich, pour tout $p\geq1$.
\endth 

{\pc DEMONSTRATION}: ${\bf \Omega}^1$ est un faisceau pour la topologie
\'etale sur Sm$(k)$. En effet, pour tout sch\'ema $X$ lisse sur $k$, $\Omega_X^1$ est un faisceau
pour la topologie de Zariski sur $X$ et il suffit de v\'erifier (cf [Mi] p50, prop 1.5) que pour
tout morphisme \'etale $X=\sp A\rightarrow Y=\sp B$ entre des sch\'emas affines sur $k$, la suite
$$\Omega^1_Y(Y)\longrightarrow\Omega^1_X(X)
\build{}_{\displaystyle\longrightarrow}^{\displaystyle\longrightarrow}
\Omega^1_{X\times_YX}(X\times_YX)$$ est exacte.

Les anneaux $A$ et $B$ sont lisses sur $k$, $A$ est \'etale sur $B$, donc
 $\Omega^1_{B/k}$ est libre sur $B$ et
$\Omega^1_{A\otimes_BA/k}=(A\otimes_B A)\otimes_B\Omega^1_{B/k}$ est libre sur $A\otimes_B A$.
L'exactitude de la suite s'en d\'eduit directement.


La topologie de Nisnevich \'etant moins fine que la topologie \'etale, ${\bf \Omega}^1$ est un
faisceau pour la topologie de Nisnevich sur Sm$(k)$.

De m\^eme, ${\bf \Omega}^p$ est un faisceau pour la topologie \'etale, par application directe du th\'eor\`eme de changement de base pour les alg\`ebres ext\'erieures (cf [Ma], Appendice C) . On obtient donc un faisceau pour la topologie de Nisnevich. 
$\diamond$

\smallskip
\th COROLLAIRE 1.2
\enonce Le complexe ${\bf \Omega}^0\rightarrow{\bf \Omega}^1\rightarrow\cdots\rightarrow{\bf \Omega}^p\rightarrow\cdots$ est un complexe de faisceaux pour la topologie de Nisnevich.
\endth 

\medskip
{\bf 2  Un th\'eor\`eme de descente}
\smallskip
Il s'agit \`a pr\'esent de munir les faisceaux ${\bf \Omega}^p$ de transferts. D'apr\`es un r\'esultat de Suslin et Voevodsky ([SV96]), tout faisceau v\'erifiant la propri\'et\'e de descente galoisienne peut \^etre muni de transferts, via la trace. Or les faisceaux ${\bf \Omega}^p$ ne v\'erifient pas cette propri\'et\'e, m\^eme si l'on se ram\`ene \`a une base lisse, comme le montre l'exemple ci-dessous (remarque 2.2) transmis par L. Illusie.

\smallskip
Il est n\'eanmoins possible d'adapter la d\'emonstration de [SV96] en g\'en\'eralisant un th\'eor\`eme de Knighten ([K73]) sur le bidual de  $ \Omega_X^p$ en tant que ${\cal O}_X$-module, consid\'er\'e ici comme un faisceau pour la topologie de Zariski sur $X$. 

Nous nous pla\c cons dans la situation  o\`u $X$ est un sch\'ema irr\'eductible et lisse sur $k$,  $W$ la normalis\'ee de $X$ dans une extension galoisienne mod\'er\'ement ramifi\'ee du corps des fonctions $K(X)$ de $X$, de groupe de
Galois $G$ et commen\c cons par construire une application $\Omega^p_W(W)\rightarrow\Omega^p_X(X)$, qui co\"\i ncide avec la trace lorsque $W$ est lisse.


\tha THEOREME  2.1
\enonce Soit $X$ est un sch\'ema irr\'eductible et normal
sur $k$ et supposons que $W$ soit la normalis\'ee de $X$ dans une
extension galoisienne mod\'er\'ement ramifi\'ee du corps des fonctions $K(X)$ de $X$, de groupe de
Galois $G$. Notons  $(\Omega^p_X)^{\ast\ast}$ (respectivement  $(\Omega^p_W)^{\ast\ast}$) le bidual
de $\Omega^p_X$ (respectivement $\Omega^p_W$) en tant que ${\cal O}_X$-module (respectivement 
${\cal O}_W$-module); alors il existe une application naturelle
$$(\Omega^p_X)^{\ast\ast}(X)\rightarrow(\Omega^p_W)^{\ast\ast}(W),$$ qui induit un isomorphisme $(\Omega^p_X)^{\ast\ast}(X)\simeq(\Omega^p_W)^{\ast\ast}(W)^G$.
\endth

\smallskip\goodbreak
{\pc REMARQUES} 2.2:  

$i)$ C. Knighten ([K73]) en a publi\'e une d\'emonstration d\'etaill\'ee pour $\Omega^1$ dans le cas
o\`u $k$  est alg\'ebriquement clos et propose \'egalement des contre-exemples dans des cas de
ramification sauvage.

$ii)$ Lorsque $X$ est lisse sur $k$, $\Omega^p_X$ est un faisceau localement libre sur $X$ et l'application naturelle $\Omega^p_X\rightarrow(\Omega^p_X)^{\ast\ast}$ est un isomorphisme pour tout entier $p\geq1$; on peut alors d\'efinir l'application suivante, que l'on note $\alpha_{W/X}$:
$$\Omega^p_W(W)\rightarrow (\Omega^p_W)^{\ast\ast}(W)\ \dfleche{{1\over\sharp
G}\sum_{g\in G}g^\ast}{}\ ((\Omega^p_W)^{\ast\ast}(W)^G\ \dfleche{\kern-2pt\displaystyle\sim}{}\ (\Omega^p_X)^{\ast\ast}(X)\ \gfleche{\kern-2pt\displaystyle\sim}{}\ \Omega^p_X(X)\quad.$$

$iii)$ {\sl Exemple:} posons $A=k[u,v]$ et $B=A[x,y]/(x^3-uv^2,y^3-u^2v,x^2-yv,y^2-xu,xy-uv)$; alors $B$ est un anneau libre sur $A$, l'extension des corps de fractions est une extension galoisienne de degr\'e $3$, de groupe de Galois $G=\bf \mu_3$ lorsque $k$ contient les racines troisi\`emes de l'unit\'e.  On peut montrer qu'il existe un \'el\'ement primitif  (plus pr\'ecis\'ement $B=A[G](1+x+y)$), mais $(\Omega^1_{B/A})^G\not=\{0\}$: en effet la forme $xdy$ est invariante par $G$.

\medskip
{\pc DEMONSTRATION }:  La d\'emonstration est calqu\'ee sur celle de C. Knighten ([K73]). 

$a)$ {\sl Construction de l'application $(\Omega^p_X)^{\ast\ast}(X)\rightarrow(\Omega^p_W)^{\ast\ast}(W)$} 

Pour $\cal M$ et $\cal N$ deux faisceaux de ${\cal O}_X$-modules, notons :

--  ${\cal H}om_X({\cal M},{\cal N})$ le faisceau $U\mapsto\ho_{{\cal O}_U}({\cal M}_{|U},{\cal N}_{|U})$,

-- ${\cal A}lt^p_X({\cal M},{\cal N})$ le faisceau associ\'e au pr\'efaisceau $U\mapsto \alt^p_{{\cal O}_U}({\cal M}_{|U},{\cal N}_{|U})$, des formes $p$-diff\'erentielles altern\'ees de $\cal M$ dans $\cal N$,

-- ${\cal D}_k({\cal O}_U,{\cal M})$ le faisceau associ\'e au pr\'efaisceau $U\mapsto D_k({\cal O}_X(U),{\cal M}_{|U})$, des $k$-d\'erivations de ${\cal O}_X$ dans $\cal M$.

Le lemme suivant n'est qu'une r\'e\'ecriture des propri\'et\'es universelles des formes diff\'erentielles de K\"ahler (cf [K73], lemme 2, p 67) et des alg\`ebres ext\'erieures.

\th LEMME 2.3
\enonce $$\eqalign{&i)\quad(\Omega^1_X)^\ast={\cal H}om_X(\Omega^1_X,{\cal O}_X)={\cal D}_k({\cal O}_X,{\cal O}_X)\cr
&ii)\quad(\Omega^p_X)^\ast={\cal H}om_X(\Omega^p_X,{\cal O}_X)={\cal A}lt^p_X(\Omega_X^1,{\cal O}_X).}$$
\endth

Puisque sur un sch\'ema normal les faisceaux duaux sont sans torsion, on peut en donner  une description explicite des \'el\'ements. Pour tout ouvert $U$ de $X$, on a
$$(\Omega^1_X)^\ast(U)=\{D\in D_k(K(X),K(X))\hbox{ tel que }\forall x\in U,  D({\cal O}_x)\subset{\cal O}_x\}$$
et similairement
$$(\Omega^p_X)^\ast(U)=\{D\in\alt^p_k(\Omega^1_{K(x)/k},K(X))\hbox{ tel que }\forall x\in U, D((\Omega^1_{X,x})^p)\subset{\cal O}_x\}.$$
Puisque $(\Omega^p_X)^{\ast\ast}(U)$ est sans torsion, l'application naturelle $(\Omega^p_X)^{\ast\ast}(U)\rightarrow \Omega^p_{K(X)/k}$ est injective et la description pr\'ec\'edente permet de reconna\^\i tre les \'el\'ements de $\Omega^p_{K(X)/k}$ provenant de $(\Omega^p_X)^{\ast\ast}(U)$. Nous pouvons \'ecrire la propri\'et\'e caract\'eristique des \'el\'ements du bidual
$$(\Omega^1_X)^{\ast\ast}(U)=\{\omega\in\Omega^1_{K(X)/k}\hbox{ tel que }\forall x\in U, \forall D\in D_k({\cal O}_x,{\cal O}_x), D(\omega)\in{\cal O}_x\}\leqno{(P_1(U))}$$
et similairement
$$(\Omega^p_X)^{\ast\ast}(U)=\{\omega\in\Omega^p_{K(X)/k}\hbox{ tel que }\forall x\in U, \forall D\in\alt^p_k(\Omega^1_{X,x},{\cal O}_x), D(\omega)\in{\cal O}_x\}.\leqno{(P_p(U))}$$
C'est la propri\'et\'e $P_1$ qui d\'efinit les diff\'erentielles de Zariski dans [K73].

Comme l'application $f:W\rightarrow X$ est dominante, elle induit une extension de corps $K(X)\rightarrow K(W)$. Pour construire l'application $(\Omega^p_X)^{\ast\ast}(X)\rightarrow(\Omega^p_W)^{\ast\ast}(W)$, il suffit de s'assurer que, pour $p\geq1$, si $\omega\in\Omega^p_{K(X)/k}$ v\'erifie la propri\'et\'e $P_p(U)$ pour un certain ouvert $U$ de $X$, son image
 $\omega'$ par l'application canonique $\Omega^p_{K(X)/k}\rightarrow \Omega^p_{K(W)/k}$ v\'erifie la propri\'et\'e $P_p(f^{-1}(U))$ sur $W$.

On utilise le fait que le bidual $(\Omega^p_X)^{\ast\ast}$ (respectivement 
$(\Omega^p_W)^{\ast\ast}$) est un faisceau de modules r\'eflexifs pour la topologie de Zariski sur
$X$ (respectivement sur $W$), c'est-\`a-dire, pour $U$ ouvert de $X$, on a
$$(\Omega^p_X)^{\ast\ast}(U)=\cap_{{{\cal P}\in U}, \dim{\cal P}=1}(\Omega^p_X)^{\ast\ast}_{\cal
P}$$ et de m\^eme pour les ouverts de $W$. Il suffit alors de montrer la propri\'et\'e en des points de
codimension $1$ et, comme $X$ et $W$ sont normaux, les anneaux locaux des points de codimension $1$ sont des anneaux de valuation discr\`ete. On se ainsi ram\`ene localement \`a $X=\sp A$ et $W=\sp B$ o\`u $A$ et $B$ sont des anneaux de valuation discr\`ete (cf [K73]).
$\diamond$
\goodbreak
\smallskip
$b)$ {\sl L'isomorphisme $(\Omega^p_X)^{\ast\ast}(X)\simeq(\Omega^p_W)^{\ast\ast}(W)^G$}

A nouveau, la d\'emonstration repose sur le fait que $(\Omega^p_X)^{\ast\ast}$ et
$(\Omega^p_W)^{\ast\ast}$ sont des faisceaux de modules r\'eflexifs et il suffit de montrer l'isomorphisme lorsque $X=\sp A$ et $W=\sp B$ o\`u $A$ et $B$ sont des anneaux de valuation discr\`ete. Dans cette situation $\Omega^p_{A/k}$ (resp  $\Omega^p_{B/k}$) est un $A$-module libre de type fini (resp un $B$-module libre) et s'identifie \`a son bidual. 
$\diamond$

\medskip
{\bf 3  Transferts}
\smallskip
Si $X_1$ et $X_2$ sont des sch\'emas lisses sur $k$ et $Z\subset X_1\times_kX_2$ un sous-sch\'ema fini sur $X_1$ et surjectif sur une composante irr\'eductible de $X_1$, il s'agit de construire, pour tout $p\geq1$, une application
$T_Z:\Omega^p_{X_2}(X_2)\rightarrow \Omega^p_{X_1}(X_1)$. L'application $Z\hookrightarrow  X_1\times_kX_2\rightarrow X_2$ compos\'ee de l'injection et de la projection sur $X_2$ d\'efinit une
fl\`eche $\Omega^p_{X_2}(X_2)\rightarrow\Omega^p_Z(Z)$. Il suffit de  construire une
application $T_{Z/{X_1}}:\Omega^p_Z(Z)\rightarrow\Omega^p_{X_1}(X_1)$ pour tout morphisme $Z\rightarrow X_1$ fini et surjectif sur une composante irr\'eductible de $X_1$ et de v\'erifier la compatibilit\'e de ces applications avec la composition des correspondances.
\smallskip
La suite de cette partie consiste \`a v\'erifier que la d\'emonstration de [SV96] s'adapte dans notre situation.

\medskip
$3.1$ {\sl Construction de $T_{Z/X}$}

Pour $X$ lisse et irr\'eductible sur $k$,  $Z\rightarrow  X$  irr\'eductible, fini et surjectif sur $X$,  on note $W$ la normalis\'ee de $X$ dans une extension
galoisienne de $K(X)$ contenant le corps des fonctions $K(Z)$ de $Z$. Alors $\ho_X(W,Z)\not=\emptyset$ et on pose (cf [SV96])
$$T_{Z/X}=\alpha_ {W/X}\circ(\sum_{q\in\ho_{X}(W,Z)}q^*)\quad,$$
o\`u $\alpha_{W/X}:\Omega^p_W(W)\rightarrow\Omega^p_{X}(X)$ est l'application d\'efinie ci-dessus ({\pc REMARQUE} 2.2.).

\th LEMME $3.1.1$ 
\enonce L'application $T_{Z/X}$ ne d\'epend pas du sch\'ema normal $W$ choisi. 
\endth

{\pc DEMONSTRATION}: consid\'erons $W_1$ et $W_2$ deux sch\'emas normaux sur
$X$ tels que les extensions $K(W_i)$ soient galoisiennes sur $K(X)$,
$K(Z)\subset K(W_i)$ pour $i=1,2$ et tels qu'il existe une application $f:
W_1\rightarrow W_2$ finie et surjective. Il s'agit de montrer que le diagramme suivant est commutatif:

$$\xymatrix{ &\Omega^p_{W_1}(W_1)\ar[dr]^{\alpha_{W_1/X}}\\
\ar @{}[r]|(0.65){(1)}\Omega^p_Z(Z)\ar[ur]^<<<<{\sum_{q\in\ho_X(W_1,Z)}q^\ast}
\ar[dr]_<<<<{\sum_{q\in\ho_{X_1}(
W_2,Z)}q^\ast\ \ }&\ar @{}[r]|(0.4){(2) }&\Omega^p_X(X)\\
&\Omega^p_{ W_2}( W_2)\ar[uu]_ {f^*}\ar[ur]_{\alpha_{W_2/X}}}$$

$(1)$ comme $Z$ est irr\'eductible et $W_i\rightarrow X$ est un recouvrement  pseudo-galoisien (au sens de [SV96]) tel que $\ho_X(W_i,Z)\not=\emptyset$, l'application naturelle 
$$\ho_X(W_i,Z)\rightarrow\ho_{K(X)}(K(Z),K(W_i))$$ est bijective et le cardinal de
$\ho_X( W_i,Z)$ est \'egal au degr\'e s\'eparable de l'extension $[K(Z):K(X)]$ ([SV96], {\sl lemme 5.7} ), d'o\`u le fait
que l'application
$$\matrix{ \ho_X(W_2,Z) &\longrightarrow &\ho_X(W_1,Z)\cr
q &\longmapsto &q\circ f}$$ est une bijection et le triangle $(1)$ commutatif.

$(2)$ Puisque $f:W_1\rightarrow W_2$ est un morphisme dominant entre sch\'emas normaux , l'application
$$f^\ast:(\Omega^p_{ W_2})^{**}(W_2)\rightarrow (\Omega^p_{W_1})^{**}(W_1)$$
existe par le th\'eor\`eme 2.1.
On note $G_1$ (respectivement $G_2$) le groupe de Galois de $K(
W_1)$ sur $K(X)$ (respectivement de $K(W_2)$ sur $K(X)$). On est dans la situation o\`u $K(W_2)\subset K(W_1)$; le groupe $G_2$ est alors un quotient de $G_1$ et l'extension $K( W_1)/K(W_2)$
est galoisienne de groupe de Galois $H=\{h\in G_1$ tel que $h(x)=x$ pour $x\in K(W_2)\}$. La fonctorialit\'e s'exprime par la
commutativit\'e du diagramme suivant:

$$\xymatrix{
& (\Omega^p_X)^{**}(X)\simeq\Omega^p_X(X) \ar[dl]_{\sim}\ar[dr]^{\sim}\\
 ((\Omega^p_{W_1})^{**}(W_1))^{G_1}  &&((\Omega^p_{ W_2})^{**}(W_2))^{G_2} \\
 (\Omega^p_{ W_1})^{**}( W_1)  \ar[u]_{{1\over\sharp G_1}\sum_{g\in G_1}g^\ast}\ar@/^/[rr]^{ {1\over\sharp H}\sum_{h\in H}h^\ast}&& 
 (\Omega^p_{W_2})^{**}( W_2)\ar[u]^{{1\over\sharp G_2}\sum_{g\in G_2}g^\ast}\ar@/^/[ll]^{f^*}\ar[r]^\sim &((\Omega^p_{W_1})^{**}(W_1))^H\\
\Omega^p_{W_1}( W_1)\ar[u] && 
\Omega^p_{W_2}( W_2)\ar[ll]_{f^*}\ar[u]}$$
\hfill $\diamond$

 \goodbreak
\smallskip
{\pc REMARQUE}: ce dernier lemme est valable aussi sur un corps de caract\'eristique $p>0$,
auquel cas, on suppose simplement que l'extension $K(W)/K(X)$ est normale de degr\'e non divisible par $p$ et on modifie l'expression de $T_{Z/X}$ par (cf [SV96])
$$T_{Z/X}=[K(W):K(X)]_{\rm ins\acute ep}\sum_{q\in\ho_X( W,Z)}q^\ast\quad.$$

\medskip
$3.2$ {\sl Extension \`a $Z$ non irr\'eductible}

Soit $Z$ un sous-sch\'ema de $X_1\times_k X_2$ fini et surjectif sur $X_1$ et notons $Z_i$ pour $1\leq i\leq d$, les composantes irr\'eductibles de $Z$ de multiplicit\'es $n_i$. Le cycle associ\'e \`a $Z$ est alors $[Z]=\sum_{1\leq i\leq d}n_i[Z_i]$ et on \'etend la d\'efinition pr\'ec\'edente par lin\'earit\'e:
$$T_{[Z]}=\sum_{1\leq i\leq d}n_iT_{Z_i}:\Omega^p_{X_2}(X_2)\ \dfleche{}{}\ \oplus_{1\leq i\leq d}\Omega^p_{Z_i}(Z_i)\ \dfleche{\sum_{1\leq i\leq d}n_iT_{Z_i/X_1}}{}\ \Omega^p_{X_1}(X_1).$$

\medskip
{\bf 4 Composition des correspondances}
\smallskip
Il reste \`a v\'erifier que le faisceau ${\bf \Omega}^p$, muni des transferts construits ci-dessus est un faisceau sur la cat\'egorie $\smc(k)$. 

Soient $X_1$, $X_2$ et $X_3$ des sch\'emas lisses et irr\'eductibles sur $k$, $Z\subset X_1\times X_2$ un sous-sch\'ema irr\'eductible, fini et surjectif sur $X_1$ et $Z'\subset X_2\times X_3$ un sous-sch\'ema irr\'eductible, fini et surjectif sur $X_2$. Par composition, on peut d\'efinir une application
$T_{Z'}\circ T_Z:\Omega^p_{X_1}(X_1)\rightarrow\Omega^p_{X_3}(X_3)$.

D'autre part, si  $p_{13}:X_1\times X_2\times X_3\rightarrow X_1\times X_3$ est la projection, on note $p_{13\ast}(
[Z\times_{X_2}Z'])=[Z']\circ [Z]$ (composition des cycles au sens de [V00]) et le morphisme $p_{13\ast}(Z\times_{X_2}Z')\rightarrow X_1$ est fini et surjectif . Ceci permet de d\'efinir l'application 
$$T_{[Z']\circ[ Z]}:\Omega^p_{X_1}(X_1)\rightarrow\Omega^p_{X_3}(X_3)\quad.$$

\tha PROPOSITION 4
\enonce Sous les notations pr\'ec\'edentes, $T_{Z'}\circ T_ Z=T_{[Z']\circ [Z]}$.
\endth
La suite de ce chapitre est consacr\'ee \`a la d\'emonstration de cette proposition.

\medskip
$4.1$ {\sl Image directe}

Soit $X$ un sch\'ema lisse et irr\'eductible sur $k$, $Z'$ irr\'eductible et $Z\dfleche{p}{}Z'\rightarrow X$ des morphismes finis et surjectifs. Notons $Z_i$ les composantes irr\'eductibles de $Z$,  $n_i$ leurs multiplicit\'es et $d_i=[K(Z_i):K(Z')]$.

\th  LEMME 4.1.1
\enonce Le diagramme suivant est commutatif 
$$\xymatrix{\oplus_{1\leq i\leq d}\Omega^p_{Z_i}(Z_i) \ar[rr]^(0.55){\sum_{1\leq i\leq d}n_iT_{Z_i/X}}&&\Omega^p_X(X)\\
&\Omega^p_{Z'}(Z')\ar[ul]^{p^\ast}\ar[ur]_{T_{p_\ast([Z])}}}$$
avec $p_\ast([Z])=\sum_{1\leq i\leq d}d_in_i[p(Z_i)]=\sum_{1\leq i\leq d}d_in_i[Z']$. 
\endth

{\pc DEMONSTRATION}:
On introduit une normalisation $W$ de $X$ dans une extension galoisienne de $K(X)$ contenant \`a la fois $K(Z')$ et tous les $K(Z_i)$ et l'on se ram\`ene \`a montrer que le triangle suivant, o\`u $n=\sum_{1\leq i\leq d}n_id_i$,
$$\xymatrix{\oplus_{1\leq i\leq d}\Omega^p_{Z_i}(Z_i) \ar[rr]^{\sum_{1\leq i\leq d\atop q_i\in\ho_X(W,Z_i)}n_iq_i^\ast} &&\Omega^p_W(W)\\
&\Omega^p_{Z'}(Z')\ar[ul]^{p^\ast}\ar[ur]_>>>>>{\ \ n\sum_{q'\in\ho_X(W,Z')}q'^\ast} }$$\nobreak
est commutatif. 

Or, puisque $W$ est normal et $W\rightarrow X$ est pseudo-galoisien, il en est de m\^eme de $W\rightarrow Z'$ et $\ho_{Z'}(W,Z_i)=\ho_{K(Z')}(K(Z_i),K(W))$ est de cardinal $d_i$ ([SV96], {\sl lemmes 5.3 et 5.7}). Tout $q'\in\ho_X(W,Z')$ fix\'e induit une application $\ho_{Z'}(W,Z_i)\rightarrow\ho_X(W,Z_i)$ injective et la compos\'ee
$$\xymatrix{&\Omega^p_{Z'}(Z')\ar[r]^{p_i^\ast}&\Omega^p_{Z_i}(Z_i)\ar[rrr]^{\sum_{q_i\in\ho_X(W,Z_i)}q_i^\ast}&&&\Omega^p_W(W)}$$ est l'application $d_iq'^\ast$,  d'o\`u le r\'esultat.$\diamond$

\medskip
{\pc CONSEQUENCE :} 
le lemme s'applique dans la situation suivante. Pour $X_1$, $X_2$ et $X_3$ des sch\'emas lisses et irr\'eductibles sur $k$, $Z\subset X_1\times X_2$ un sous-sch\'ema irr\'eductible, fini et surjectif sur $X_1$ et $Z'\subset X_2\times X_3$ un sous-sch\'ema irr\'eductible, fini et surjectif sur $X_2$, 
notons $Z_i$,  pour $1\leq i\leq d$, les composantes irr\'eductibles de $Z\times_{X_2}Z'$ de multiciplit\'es $n_i$ et de m\^eme $Z'_j$, pour $1\leq j\leq f$, les composantes irr\'eductibles de $p_{13}(Z\times_{X_2}Z')$ de multiciplicit\'es $n'_j$. On a alors
 $[Z\times_{X_2}Z']=\sum_{1\leq i\leq d}n_i[Z_i]$, 
$$p_{13\ast}([Z\times_{X_2}Z'])=\sum_{1\leq i\leq d}n_ip_{13\ast}[Z_i]=\sum_{1\leq j\leq f}n'_j[Z'_j]$$ 
et la commutativit\'e du diagramme
$$\xymatrix{\oplus_{1\leq i\leq d}\Omega^p_{Z_i}(Z_i) \ar[rr]^{\sum_{1\leq i\leq d}n_iT_{Z_i/X_1}}{}&&\Omega^p_{X_1}(X_1)\\
&\oplus_{1\leq j\leq f}\Omega^p_{Z'_j}(Z'_j)\ar[ul]^{ p_{13}^\ast}\ar[ur]_{\ \ \ \ \ \sum_{1\leq j\leq f}n'_jT_{Z'_j/X_1}\quad.}}$$

\medskip
$4.2$ {\sl Composition}
\th PROPOSITION 4.2.1
\enonce Introduisons $W'$ (resp $W$) une normalisation de $X_2$ (resp de $X_1$ ) dans une extension galoisienne de $K(X_2)$ (resp de $K(X_1)$) contenant $K(Z')$ (resp $K(Z)$ et $K(Z_i)$ pour tout $i$, $1\leq i\leq d$), de groupe de Galois $G'=\Gal(K(W')/K(X_2))$ (resp $G=\Gal(K(W)/K(X_1))$). Le diagramme suivant est commutatif:

$$\xymatrix{&&&\oplus_{1\leq i\leq d}\Omega^p_{Z_i}(Z_i) \ar[drrr]^(0.6){\sum_{1\leq i\leq d\atop q_i\in\ho_{X_1}(W,Z_i)}n_iq_i^\ast}\\
\Omega^p_{Z'}(Z')\ar[urrr]^{p_{23}^\ast}\ar[dr]^(0.65){\ \sum_{q\in\ho_{X_2}(W',Z')}q^\ast} 
&&&&&&\Omega^p_W(W)\\
&\Omega^p_{W'}(W')\ar[r]^(0.45){\alpha} 
&(\Omega^p_{W'})^{\ast\ast}(W')\ar[r]^{{1\over\sharp G'}\sum_{g\in G'}g^\ast}
&(\Omega^p_{W'})^{\ast\ast}(W')^{G'}\ar[r]^(0.55)\sim 
&\Omega^p_{X_2}(X_2)\ar[r]^{p_2^\ast}
&\Omega^p_Z(Z)\ar[ur]^(0.4){\sum_{q\in\ho_{X_1}(W,Z)}q^\ast}}
$$ 
\endth

{\pc DEMONSTRATION}: la d\'emonstration suit pas-\`a-pas celle de [SV96]; il s'agit juste de v\'erifier qu'elle reste valable malgr\'e le passage par le bidual.

Choisissons $W'_0$ une composante irr\'eductible de $W'\times_{X_2}Z$ et prenons $W$ la normalisation de $W'_0$ dans une extension galoisienne de $K(X_1)$ contenant $K(W'_0)$ et $K(Z_i)$ pour $1\leq i\leq d$. Nous sommes dans la situation suivante, o\`u les fl\`eches verticales sont finies et surjectives sur une composante irr\'eductible du but.

$$\xymatrix{ W\ar[d]^{s}\ar@/_1.5pc/[dd]_{q_i}\\ 
W'_0 \ar @<-0.3ex>@{^{(}->}[r]\ar[d]^{r'_i}\ar@/^1.5pc/[rr]^{\pi'} &W'\times_{X_2}Z \ar[r]\ar[d]^{r'} 
&W' \ar[d]^{q'}\ar@/^1.5pc/[dd]^{p'}\\
Z_i \ar[dr]\ar@<-0.3ex> @{^{(}->}[r]&Z'\times_{X_2}Z \ar[r]^(0.6){p_{23}}\ar[d]^{p_{12}} 
&Z'\ar[d]^{p'_1}\\
&Z \ar[r]^{\pi}\ar[d]^{p_1} &X_2\\
&X_1}
$$

De plus,

-- les fl\`eches $\pi':W'_0\rightarrow W'$ et $p'=p'_1\circ q':W'\rightarrow X_2$ ne d\'ependent pas du choix de $q'\in\ho_{X_2}(W',Z')$;

-- pour chaque $q'\in\ho_{X_2}(W',Z')$, il existe $r':W'\times_{X_2}Z\rightarrow Z'\times_{X_2}Z $ unique par propri\'et\'e universelle du produit fibr\'e;

-- comme $W'_0$ est irr\'eductible, pour chaque $r'$,  il existe $i$, $1\leq i\leq d$ unique tel que $r'(W'_0)=Z_i=r'_i(W)$ et $r'_i$ est fini et surjectif.

\th LEMME 4.2.2
\enonce
Pour $a\in\Omega^p_{Z'}(Z')$, 
$$\sum_{q\in\ho_{X_1}(W,Z)}q^\ast\circ p_2^\ast\circ T_{Z'/X_2}(a)={1\over [K(W'_0):K(Z)]}\sum_{s\in\ho_{X_1}(W,W'_0)}s^\ast\circ\pi'{}^\ast\circ p'{}^\ast\circ T_{Z'/X_2}(a).$$
\endth
{\pc DEMONSTRATION}:
$$\sum_{q\in\ho_{X_1}(W,Z)}q^\ast\circ p_2^\ast\circ T_{Z'/X_2}(a)
=\sum_{q\in\ho_{X_1}(W,Z)}s^\ast\circ\pi'{}^\ast\circ p'{}^\ast\circ T_{Z'/X_2}(a).$$
Par le lemme 5.9 de [SV96], on sait que $W'_0/Z$ est pseudo-galoisien de groupe 
$$G'=\ho_{K(Z)}(K(W'_0),K(W'_0))=\ho_{K(Z)}(K(W'_0),K(W))\quad;$$  de plus, comme $W$ est normal
$$\ho_{X_1}(W,Z)=\ho_{K(X_1)}(K(Z),K(W))$$
et
$$\ho_{X_1}(W,W'_0)=\ho_{K(X_1)}(K(W'_0),K(W)).$$
L'application
$$\matrix{ \ho_{K(X_1)}(K(W'_0),K(W)) &\longrightarrow &\ho_{K(X_1)}(K(Z),K(W))\cr
s &\longmapsto &q=s_{|K(Z)}}$$ 
est surjective et ses fibres sont de cardinal $[K(W'_0):K(Z)]$.
En effet, $s$ et $s'$ ont m\^eme image si et seulement s'il existe $t\in\ho_{K(Z)}(K(W'_0),K(W'_0))$ tel que $s'=s\circ t$  et $\ho_{K(Z)}(K(W'_0),K(W'_0))$ est de cardinal $[K(W'_0):K(Z)]$.

On en d\'eduit le lemme.$\diamond$

\smallskip
\th LEMME 4.2.3
\enonce $$p'{}^\ast\circ T_{Z'/X_2}(a)=\sum_{q'\in\ho_{X_2}(W',Z')}q'{}^\ast(a)\quad.$$
\endth
{\pc DEMONSTRATION}: 
Rappelons que $p'{}^\ast\circ T_{Z'/X_2}(a)=p'{}\ast\circ({1\over\sharp G'}\sum_{g\in G'}g^\ast)\circ\alpha\circ\sum_{q\in\ho_{X_2}(W',Z')}q^\ast$. Par fonctorialit\'e $g^\ast\circ\alpha=\alpha\circ g^\ast$ pour $g\in G'=\ho_{X_2}(W',W')$; d'autre part, le noyau de $\alpha$ est $(\Omega^p_{W'})_{\rm tors}$ et, comme $X_2$ est lisse, $\Omega^p_{X_2}(X_2)$ n'a pas de torsion, on en d\'eduit que la compos\'ee 
$$\xymatrix{\Omega^p_{W'}(W') \ar[r]^\alpha &(\Omega^p_{W'})^{\ast\ast}(W') \ar[rr]^{{1\over\sharp G'}\sum_{g\in G'}g^\ast} &&((\Omega^p_{W'})^{\ast\ast}(W'))^{G'}\ar[d]^\sim\\
&&&\Omega^p_{X_2}(X_2)\ar[ulll]_{p'{}^\ast}}$$
induit l'identit\'e sur $\Omega^p_{X_2}(X_2)$, d'o\`u le r\'esultat.$\diamond$
\smallskip
On obtient alors
$$\sum_{q\in\ho_{X_1}(W,Z)}q^\ast\circ p_2^\ast\circ T_{Z'/X_2}(a)={1\over [K(W'_0):K(Z)]}\sum_{s\in\ho_{X_1}(W,W'_0)\atop q'\in\ho_ {X_2}(W',Z')}s^\ast\circ\pi'{}^\ast\circ q'{}^\ast(a).$$

\smallskip
La fin du calcul repose sur le r\'esultat suivant:
\th LEMME 4.2.4
\enonce Pour $i$, $1\leq i\leq d$, notons $a_i$ le pull-back de $a$ dans $\Omega^p_{Z_i}(Z_i)$, de telle sorte que $\pi'{}^\ast\circ q'{}^\ast(a)=\sum_{1\leq i\leq d}r'_i{}^\ast(a_i)$. Alors
$$\sum_{s\in\ho_{X_1}(W,W'_0)\atop q'\in\ho_ {X_2}(W',Z')}s^\ast\circ\pi'{}^\ast\circ q'{}^\ast(a)=\sum_{1\leq i\leq d \atop {r'_i\in\ho_ Z(W'_0,Z_i)\atop q_i\in\ho_{X_1}(W,Z_i)}}n_i[K(W'_0):K(Z_i)]q_i{}^\ast(a_i).$$
\endth

{\pc DEMONSTRATION}:

Dans notre situation, les r\'esultats de [SV96] (§ 5) s'appliquent, d'o\`u

{\parindent=0.5cm
{$i)$} $n_i=\{q'\in\ho_ {X_2}(W',Z') \hbox{ tel que }q'\mapsto r'_i\}$, 

{$ii)$} $[K(Z']:K(X_2)]=\sum_{1\leq i\leq d}n_i[K(Z_i):K(Z)]$ et 

{$iii)$} $\ho_Z(W'_0,Z_i)=\ho_Z(W,Z_i)$. 
\par}

On peut transformer l'expression
$$\sum_{s\in\ho_{X_1}(W,W'_0)\atop q'\in\ho_ {X_2}(W',Z')}s^\ast\circ\pi'{}^\ast\circ q'{}^\ast(a)=\sum_{s\in\ho_{X_1}(W,W'_0)}\sum_ {1\leq i\leq d\atop r'_i\in\ho_ Z(W'_0,Z_i)}n_is^\ast\circ r'_i{}^\ast(a_i).$$
Chaque $r'_i\in\ho_Z(W'_0,Z_i)$ d\'efinit un plongement $r'_i:K(Z_i)\llhook K(W'_0)$ et une application
$$\matrix{ \ho_{K(X_1)}(K(W'_0),K(W)) &\longrightarrow &\ho_{K(X_1)}(K(Z_i),K(W))\cr
s &\longmapsto &s_{|K(Z_i)}=s\circ r'_i=q_i\quad,}$$
dont les fibres sont de cardinal $[K(W'_0):K(Z_i)]$, d'o\`u le lemme.$\diamond$

\smallskip
Pour achever la d\'emonstration de la proposition et conclure, il suffit de remarquer
que $[K(Z_i):K(Z)]$ et $\ho_Z(W,Z_i)$ ont m\^eme cardinal.

\smallskip
{\pc CONSEQUENCES}: $i)$ La construction des transferts ci-dessus commute avec la diff\'erentielle et on obtient un complexe de faisceaux de Nisnevich avec transferts, objet de $D(\sh_{Nis}(\smc(k)))$ , not\'e ${\bf\Omega}^\bullet$ et appel\'e {\sl complexe motivique de De Rham}. 

$ii)$ De plus pour tout $n$ entier, consid\'erons $$\tau_{\leq n}{\bf\Omega}^\bullet:{\bf \Omega}^1\rightarrow\cdots\rightarrow{\bf \Omega}^{n-1}\rightarrow\ke d,$$ le complexe tronqu\'e. C'est un complexe born\'e sup\'erieurement de faisceaux de Nisnevich avec transferts, dont la cohomologie est invariante par homotopie, ce qui en fait un objet de $\dme(k)$.

\bigskip
{\bf 5 R\'ealisation de De Rham}

\smallskip  
{\pc DEFINITION} 5.1   La {\sl r\'ealisation de De Rham} d'un motif effectif  ${\bf M}$ objet 
de $\dme(k)$ est le $k$-espace 
vectoriel gradu\'e dont le terme de degr\'e $q$ est
$${\bf H}_{DR}^q ({\bf M}) = \limind_n { \ho}_{\dme(k)} ({\bf M},\tau_{\leq n} \om^\bullet[q]) .$$

\smallskip
{\pc REMARQUES}  5.2: 

$i)$ Le foncteur $H^\bullet_{DR}$ ainsi construit, limite inductive filtrante de foncteurs
Hom, est un foncteur homologique. 

$ii)$ Si  ${\bf M}$ est le motif ${\bf M}(X)$ d'un sch\'ema $X$ lisse sur $k$, de dimension inf\'erieure ou \'egale \`a $n$, alors
$${\bf H}_{DR}^q ({\bf M}(X))=\ho_{\dme(k)}({\bf M}(X),\tau_{\leq n}\om^\bullet[q])$$
et
$$\ho_{\dme(k)}({\bf M}(X),\tau_{\leq n}\om^\bullet[q])={\bf
H}^q_{Nis}(X,\tau_{\leq n}\om^\bullet)={\bf H}^q_{Zar}(X,\om^\bullet)={\bf H}^q_{DR}(X)$$
( cf [V00], propositions 3.1.9 et 3.1.12)
 et l'on voit ainsi que
${\bf H}^q_{DR} ({\bf M}(X) )$ est le $q$-i\`eme groupe de cohomologie de De Rham de $X$.

\smallskip
\th PROPOSITION 5.2 
\enonce
$a)$ La r\'ealisation de De Rham du motif de Tate est 
$H^0_{DR} ({\bf Z} (1)) = k$ et $H^p_{DR} ({\bf Z} (1)) = 0$ pour $p\not=0$.

 $b)$ La r\'ealisation de De Rham induit un foncteur de 
la cat\'egorie $\dm (k)$ dans la cat\'egorie des $k$-espaces vectoriels gradu\'es.
\endth
{\pc DEMONSTRATION}

$a)$ Cela d\'ecoule directement du triangle distingu\'e
$({\bf Z} (1))[2] \rightarrow {\bf M}( {\bf P}^1) \rightarrow {\bf Z}
\rightarrow {\bf Z}(1)[1]$.

$b)$ Par le th\'eor\`eme de simplification de Voevodsky ( "cancellation theorem" cf [MVW]16.25), on a des isomorphismes $\ho_{\dme(k)}({\bf M}(p),(\tau_{\leq n}\om^\bullet[q])\otimes{\bf Z}(p))=\ho_{\dme(k)}({\bf M},\tau_{\leq n}\om^\bullet[q])$ pour un objet $M$ de $\dme(k)$ et tout entier $p\geq 1$, ce qui permet d'\'etendre le foncteur aux complexes motiviques non effectifs en posant 
$${\bf H}_{DR}^q ({\bf M}) = \limind_n { \ho}_{\dme(k)} ({\bf M}(p),\tau_{\leq n} \om^\bullet[q]\otimes{\bf Z}(p)) ,$$
pour $p$ assez grand. $\diamond$

\smallskip
\th PROPOSITION 5.3
\enonce
La r\'ealisation de De Rham pr\'esent\'ee ici co\"\i ncide avec celle d\'efinie par A. Huber sur la cat\'egorie des motifs g\'eom\'etriques.
\endth

{\pc DEMONSTRATION}: c'est une cons\'equence de la proposition 2.1.2 de [H00] et du fait que, pour un  sch\'ema projectif lisse, la r\'ealisation de De Rham du motif associ\'e  co\"\i ncide avec la composante de De Rham du foncteur de [H00]. $\diamond$

\bigskip
\biblio {Bibliographie}

 \ref{H00} \AUTHOR = {Huber, Annette}\TITLE = {Realization of {V}oevodsky's motives} \JOURNAL = {Journal of Algebraic Geometry}\VOLUME = {9}\YEAR = {2000} \NUMBER = {4}\PAGES = {755--799}\article

\ref{H95} \AUTHOR = {Huber, Annette} \TITLE = {Mixed motives and their realization in derived categories} \SERIES = {Lecture Notes in Mathematics} \VOLUME = {1604}\PUBLISHER = {Springer-Verlag}\ADDRESS = {Berlin}\YEAR = {1995}\book
		
\ref{K73}\AUTHOR = {Knighten, Carol M}\TITLE = {Differentials on quotients of algebraic varieties}\JOURNAL = {Transactions of the American Mathematical Society}\VOLUME = {177}\YEAR = {1973}\PAGES = {65--89}\articles

\ref{LW}\AUTHOR={Lecomte, Florence - Wach, Nathalie}\TITLE ={R\'ealisations des complexes motiviques de Voevodsky}\prep	

\ref{Ma}\AUTHOR = {Matsumura, Hideyuki}\TITLE = {Commutative ring theory}\SERIES = {Cambridge Studies in Advanced Mathematics}\VOLUME = {8}\PUBLISHER = {Cambridge University Press}\ADDRESS = {Cambridge}\YEAR = {1989}\book

\ref{Mi}\AUTHOR = {Milne, James S.}\TITLE = {\'{E}tale cohomology}\SERIES = {Princeton Mathematical Series}\VOLUME = {33}\PUBLISHER = {Princeton University Press}\ADDRESS = {Princeton, N.J.}\YEAR = {1980}\book

\ref{MVW}\AUTHOR = {Mazza, Carlo - Voevodsky, Vladimir - Weibel, Charles}\TITLE = {Lecture notes on motivic cohomology}\SERIES = {Clay Mathematics Monographs}\VOLUME = {2}\PUBLISHER = {American Mathematical Society}\ADDRESS = {Providence, RI}\YEAR = {2006}\book

\ref{N89}\AUTHOR = {Nisnevich, Ye. A.}\TITLE = {The completely decomposed topology on schemes and associated descent spectral sequences in algebraic {$K$}-theory}\BOOKTITLE = {Algebraic $K$-theory: connections with geometry and topology (Lake Louise, AB, 1987)}\SERIES = {NATO Adv. Sci. Inst. Ser. C Math. Phys. Sci.}\VOLUME = {279}\PAGES = {241--342}\PUBLISHER = {Kluwer Acad. Publ.}\ADDRESS = {Dordrecht}\YEAR = {1989}\incollection

\ref{SV96}\AUTHOR = {Suslin, Andrei - Voevodsky, Vladimir}\TITLE = {Singular homology of abstract algebraic varieties}\JOURNAL = {Inventiones Mathematicae}\VOLUME = {123}\YEAR = {1996}\NUMBER = {1}\PAGES = {61--94}\article

\ref{V00}\AUTHOR = {Voevodsky, Vladimir}\TITLE = {Triangulated categories of motives over a field}\BOOKTITLE = {Cycles, transfers, and motivic homology theories}\SERIES = {Ann. of Math. Stud.}\VOLUME = {143}\PAGES = {188--238}\PUBLISHER = {Princeton Univ. Press}\ADDRESS = {Princeton, NJ}\YEAR = {2000}\incollection

\end

\item{[MVW]} C. Mazza, V. Voevodsky et C.Weibel, {\it Notes on Motivic 
Cohomology}, 
preprint 2002.

\item{[SGA 4]}
 
\item{[Ver77]} J.L. Verdier, {\it Cat\'egories d\'eriv\'ees, \'etat 0} in SGA4 1/2, Lect. 
Notes in Math. 569, pp 262-311, Springer (1977).

\item{[V98]} V. Voevodsky, {\it The ${\bf A}^1$-homotopy theory}, in Proceedings of 
the 
International Congress of Mathematicians, Berlin 1998, vol. 1, pp 579-604.